\documentstyle{aipproc}
\begin{document}
\newtheorem{definition}{\normalsize\sc Definition}
\newtheorem{prop}[definition]{\normalsize\sc Proposition}
\newtheorem{lem}[definition]{\normalsize\sc Lemma}
\newtheorem{corollary}[definition]{\normalsize\sc Corollary}
\newtheorem{theorem}[definition]{\normalsize\sc Theorem}
\newtheorem{example}[definition]{\normalsize\sc  Example}
\newtheorem{remark}[definition]{\normalsize\sc Remark}

\def\sw#1{{\sb{(#1)}}}
\def\eps{{\epsilon}}
\def\id{{\rm id}}
\def\note#1{}

\title{Line bundles on quantum spheres}

\author{Tomasz Brzezi\'nski$^*$\thanks{Lloyd's Tercentenary Fellow.
On~leave from:
Department of Theoretical Physics, University of \L\'od\'z,
Pomorska 149/153, 90--236 \L\'od\'z, Poland.
{\sc E-mail}: {\tt tb10@york.ac.uk}}
and Shahn Majid$^\dagger$\thanks{{\sc E-mail}:
{\tt majid@damtp.cam.ac.uk}}}
\address{$^*$Department of Mathematics, University
of York, Heslington, York YO10 5DD, England\\
$^\dagger$ DAMTP, University of Cambridge, Cambridge CB3 9EW, England}

%\lefthead{LEFT head}
%\righthead{RIGHT head}
\maketitle

\begin{abstract}
The (left coalgebra) line bundle associated to the quantum Hopf
fibration
of any quantum two-sphere is shown to be a
 finitely generated projective
module. The corresponding projector is constructed and its monopole
charge is computed. It
is shown that the Dirac $q$-monopole connection on any quantum
two-sphere induces the Grassmannian
connection built with this projector.
\end{abstract}

\section*{Introduction}
In the standard approach to non-commutative geometry \cite{Con:non} (see
\cite{Lan:int} for an accessible introduction) the
notion of a vector bundle is identified with that of a finitely
generated projective module $E$ of an algebra $B$. Recall that $E$ is
a finitely generated projective  left
$B$-module if $E$ is
isomorphic to $B^ne := \{ (b_1,\ldots,
b_n)e\; |\;b_1,\ldots,b_n\in B\}$, where $e$ is
an $n\times n$ matrix with entries from $B$ such that $e^2 =e$ ($e$ is
called a {\em projector} of $E$). The algebra $B$ plays the
role of an algebra of functions on a manifold and $E$ is thought of as a
module of sections (matter fields) on a vector bundle over this
manifold. A (universal) connection or a covariant derivative on $E$
is a map $\nabla : E\to \Omega^{1}B\otimes_B E$,
where $\Omega^1B$ denotes the space of (universal) differential
$1$-forms on $B$,
such that for all $\rho \in \Omega^1B$, $x\in E$,
\begin{equation}
\nabla(\rho \otimes x) = \rho\nabla(x) +d\rho\otimes x.
\label{nabla}
\end{equation}

A different approach to non-commutative gauge theories, based on
quantum groups,  was introduced in \cite{BrzMaj:gau}.
In this approach one begins with a quantum principal bundle $P$
constructed algebraically as a Hopf-Galois extension of $B$ by a
Hopf algebra (quantum group) $H$ (see \cite{Mon:hop} for a review
of Hopf-Galois extensions). One then constructs a vector bundle as
a module $E$ associated to $B$ and an $H$-comodule
(corepresentation) $V$. Any strong connection (gauge field) in $P$
induces then a covariant derivative in $E$. An example of a quantum
principal bundle is the quantum Hopf fibration over the standard
quantum two-sphere in \cite{BrzMaj:gau}. The Dirac
$q$-monopole is the canonical connection on this bundle. In a recent
paper \cite{HajMaj:pro} the line bundle associated to the
$q$-monopole bundle has been constructed. It has been shown that it
is a finitely generated projective module. The  covariant
derivative induced by the $q$-monopole connection is the {\em
Grassmannian connection} (see eq.~(\ref{gras}) below).

On the other hand, the standard quantum two-sphere is
only one of the infinite family of quantum two-spheres constructed
in \cite{Pod:sph}. To describe gauge theory on all such spheres one
needs the notions of a {\em coalgebra principal bundle}, introduced
in \cite{BrzMaj:coa}, and an {\em associated module} as an
algebraic version of a coalgebra vector bundle, introduced in
\cite{Brz:mod}. In this note we announce some of the results of a
forthcoming paper \cite{BrzMaj:geo}, in which the theory of
\cite{BrzMaj:coa} is generalised and combined with \cite{Brz:mod}
to give the theory of connections on coalgebra principal and
associated bundles, illustrated by the example of monopole
connection for all quantum two-spheres. The new results of this
note then consist of an explicit projective module
description of line bundles over all quantum two-spheres,
following the method used for the standard quantum sphere in
\cite{HajMaj:pro}.

We work over the field $k$ of real or complex numbers, and we use
the standard coalgebra notation. Thus, for a coalgebra $C$,
$\Delta$ is the coproduct and $\eps$ is the counit. In a Hopf
algebra $S$ denotes the antipode. We use Sweedler's notation
for the coproduct, $\Delta(c) = c\sw 1\otimes c\sw 2$
(summation understood). If $P$ is a right $C$-comodule then
$\Delta_P :P\to P\otimes C$ denotes the right coaction. On elements
we write $\Delta_P(p) = p\sw 0\otimes p\sw 1$ (summation
understood). If $V$ is a left $C$-comodule then ${}_V\Delta:V\to
C\otimes V$ denotes the left coaction. The reader not familiar with
coalgebras is referred to \cite{Swe:hop} or to any modern textbook
on quantum groups.

\section*{Coalgebra gauge theory}
We begin with the following definition from \cite{BrzHaj:coa} (which
generalises the earlier definition of a coalgebra principal bundle in
\cite{BrzMaj:coa}).
\begin{definition}
Let $C$ be a coalgebra,
$P$ an algebra and  a right $C$-comodule, and $B$ a subalgebra of $P$,
$B :=\{b\in P\,|\;\forall p\in P \; \Delta_P(bp)=b\Delta_P(p)\}$.
We say that
P is a {\em coalgebra-Galois extension}
(or a {\em coalgebra principal bundle}) of $B$ {\em iff} the canonical left
$P$-module  right $C$-comodule map
$can: P\otimes\sb B
P\to P\otimes
C$, $can(p\otimes p')= pp'\sw 0\otimes p'\sw 1 $
is bijective.
Such a coalgebra-Galois extension is denoted by $P(B)^C$.
\label{cge}
\end{definition}
The definition of a coalgebra principal bundle implies \cite{BrzHaj:coa}
that there exists
a unique {\em entwining map} $\psi :C\otimes P\to P\otimes C$ (cf.
\cite[Definition~2.1]{BrzMaj:coa} for a definition of an entwining map).
Explicitly,
$\psi(c\otimes p) = can\circ (can^{-1}(1\otimes c)p)$.
The properties of $\psi$ \cite[Proposition~2.2]{BrzMaj:coa}
imply that the map
$\psi^2 = (\id_P\otimes \psi)\circ(\psi\otimes\id_P): C\otimes P\otimes
P\to P\otimes P\otimes C$, restricts to $\psi^2 :C\otimes \Omega^1P\to
\Omega^1 P\otimes C$, where $\Omega^1P = \{\sum_i p_i\otimes p'_i \in
P\otimes P \; |\; \sum_ip_ip'_i =0\}$ is the bimodule of
universal one-forms on $P$. Also $\psi^2\circ(\id_C\otimes d) = (d\otimes
\id_C)\circ\psi$, where $d :P\to \Omega^1P$, $d(p) = 1\otimes p-
p\otimes 1$ is the universal differential. In other words, the
universal differential calculus on $P$ is covariant with respect to
$\psi$. This allows one to introduce the following
definition in \cite{BrzMaj:geo}, which  extends the earlier
one in \cite{BrzMaj:coa}.
\begin{definition}
Let $P(B)^C$ be a coalgebra principal bundle with the entwining map
$\psi$.
A {\em connection one-form} in $P(B)^C$ is a
linear map $\omega :C\to \Omega^1P$ such that:

(i) $1\sw 0\omega(1\sw 1) = 0$, where $1\sw 0\otimes 1\sw 1 =
\Delta_P(1)$.

(ii) For all $c\in C$, $\chi\circ\omega(c) = 1\otimes c -
\eps(c)\Delta_P(1)$, where $\chi: P\otimes P\to P\otimes C$,
$\chi(p\otimes p') = p\Delta_P(p')$.

(iii) For all $c\in C$, $\psi^2(c\sw 1\otimes \omega(c\sw 2)) =
\omega(c\sw 1)\otimes c\sw 2$.
\label{con}
\end{definition}
Conditions (i)-(iii) have the following geometric meaning. By (i),
$\omega$ can be thought of as acting on $\ker\eps$, which can be
interpreted as a dual of the ``Lie algebra" of $C$. Condition (ii)
means that an application of $\omega$ to a vector field obtained by
lifting an element of a Lie algebra gives back this element of a
Lie algebra. Finally, (iii) is the covariance of $\omega$ under the
adjoint coaction. As shown in
\cite{BrzMaj:coa}\cite{BrzMaj:geo}, connection 1-forms are
in bijective correspondence with {\em connections}, i.e., projections
$\Pi$ in $\Omega^1P$ with the fixed kernel
$P(\Omega^1B)P$ and such that $\Pi\circ d$ is right $C$-covariant. 
A connection is said to be a {\em strong
connection}, if for all $p\in P$, $dp - p\sw 0\omega(p\sw 1)\in
(\Omega^1 B)P$ \cite{Haj:str}.

Given a coalgebra principal bundle $P(B)^C$, and a left $C$-comodule
 $V$ with the coaction ${}_V\Delta :V\to
C\otimes V$, one constructs the left $B$-module
$$
E := \{\sum_i p^i\otimes v^i \in P\otimes V\; | \; \sum_i
\Delta_P(p_i)\otimes v^i = \sum_i p_i\otimes{}_V\Delta(v^i)\}.
$$
The module $E$, introduced in \cite{Brz:mod}, plays the role of a
fibre bundle associated to $P$ with fibre $V$. The coalgebra
gauge theory in \cite{BrzMaj:geo} ensures that any strong
connection $\omega$ in $P$ induces a covariant derivative
$\nabla
:E\to
\Omega^1B\otimes\sb B E$ on $E$ (in the sense of equation (\ref{nabla}))
 given by
\begin{equation}
\nabla(\sum_i p^i\otimes v^i) = \sum_i dp^i\otimes v^i -\sum_i p^i\sw
0\omega(p^i\sw 1)\otimes v^i.
\label{cov.der}
\end{equation}
Since the existence of a connection in  a module $E$ is
equivalent to $E$ being projective
\cite[Corollary~8.2]{CunQui:alg}, we conclude that if there is a
strong connection in $P(B)^C$, then every associated left
$B$-module $E$ is projective.

\section*{Left coalgebra line bundles on quantum spheres}
Recall that the quantum group $SU_q(2)$ is a free algebra
generated by $1$ and  the matrix of generators ${\bf t} =
\pmatrix{\alpha &\beta
\cr \gamma &\delta}$, subject to the relations:
$$
\alpha\beta =q\beta\alpha, \quad \alpha\gamma = q\gamma\alpha, \quad
\beta\gamma = \gamma\beta, \quad \beta\delta = q\delta\beta,
$$
$$
\gamma\delta = q\delta\gamma, \quad \alpha\delta = \delta\alpha +
(q-q^{-1})\beta\gamma, \quad \alpha\delta - q\beta\gamma = 1.
$$
$SU_q(2)$ is a Hopf algebra of a matrix type, i.e.
$$
\Delta (t_{ik}) =\sum_{j=1}^2t_{ij}\otimes t_{jk},
\quad \eps({\bf t}) = 1, \quad
S\pmatrix{\alpha &\beta \cr \gamma
&\delta} = \pmatrix{\delta &-q^{-1}\beta \cr -q\gamma &\alpha} .
$$
The quantum two-spheres $S^2_{qs}$ \cite{Pod:sph} are homogeneous
spaces of $SU_q(2)$. For a given $q$ they can be viewed as
subalgebras of $SU_q(2)$ generated by $1$ and
$$
\xi = s(\alpha^2 - q^{-1}\beta^2) +(s^2-1)q^{-1}\alpha\beta, \quad \eta =
s(q\gamma^2 - \delta^2)
+(s^2-1)\gamma\delta,
$$
$$ \zeta = s(q\alpha\gamma - \beta\delta)
+ (s^2-1)q\beta\gamma ,
$$
where $s\in [0,1]$. The standard quantum sphere corresponds to
$s=0$\footnote{It is perhaps
more customary to parametrise quantum
spheres by the parameter $c\geq 0$ \cite{Pod:sph}. The parameter
$s$ which is used here is related to $p$ in \cite{Brz:hom} via
$s-s^{-1} = p^{-1}$, and $p^2$ should be identified with $c$ 
in \cite{Pod:sph}.}.

It is shown in \cite{Brz:hom} that for any quantum
sphere there is a coalgebra principal bundle with the base $B=S^2_{qs}$
and the total space $P=SU_q(2)$. The structure
coalgebra is  a quotient
$C_s=SU_q(2)/J_s$ where
$ J_s=\{\xi-s,\eta+s,\zeta\}SU_q(2)$.
We denote by $\pi_s$ the canonical projection $SU_q(2)\to C_s$. The
coproduct on $C_s$ is obtained by projecting down the coproduct in
$SU_q(2)$ by $\pi_s$. The coaction of $C_s$ on $SU_q(2)$ is given  by
$\Delta_{SU_q(2)} = (\id_{SU_q(2)}\otimes \pi_s)\circ \Delta$.
 It has been recently found in \cite{MulSch:hom}, that the
coalgebras $C_s$ are spanned by group-like elements
($c\in C_s$ is group-like if  $\Delta(c) = c\otimes c$). These
are computed explicitly in \cite{BrzMaj:geo} and are given by
$$
\pi_s(1), \quad g^+_n=\pi_s(\prod_{k=0}^{n-1}
(\alpha + q^ks\beta)), \quad g^-_n=\pi_s(\prod_{k=0}^{n-1}
(\delta - q^{-k}s\gamma)), \quad n=1,2,\ldots
$$
(all products
increase from left to right).

Similarly to
\cite{HajMaj:pro}, consider a
one-dimensional left corepresentation $V=k$ of $C_s$ given by the coaction
${}_k\Delta : k\to C_s\otimes k$, ${}_k\Delta(1) = g^+_1 \otimes 1$.
The
left $S_{qs}^2$-module $E_s$ associated to $SU_q(2)$ and this
corepresentation comes out as
$$
E_s =\{x(\alpha+s\beta) +y(\gamma +s\delta) \; |\; x,y\in
S_{qs}^2\}.
$$
The module $E_s$ is the {\em line bundle} over $S^2_{qs}$.
Our first goal is to  construct the projector for $E_s$.

Consider the following matrix with entries from $S_{qs}^2$,
$$
e_s= {1\over 1+s^2}\pmatrix{1-\zeta & \xi \cr -\eta & s^2 +q^{-2}\zeta}.
$$
Notice that $e_s$ can be also written in the following form
\begin{eqnarray*}
e_s & = &{1\over 1+s^2}\pmatrix{(\alpha+s\beta)(\delta -qs\gamma) &
(\alpha+s\beta)(s\alpha -q^{-1}\beta) \cr (\gamma+s\delta)(\delta -
qs\gamma) &
(\gamma+s\delta)(s\alpha - q^{-1}\beta)}\\
& = &{1\over 1+s^2}
\pmatrix{\alpha+s\beta
\cr \gamma+s\delta}(\delta - qs\gamma, s\alpha-q^{-1}\beta).
\end{eqnarray*}
The use of the relation
\begin{equation}
(\delta -qs\gamma)(\alpha+s\beta) +
(s\alpha-q^{-1}\beta)(\gamma+s\delta) = 1+s^2,
\label{fact}
\end{equation}
makes it clear that $e_s^2 =e_s$, i.e., $e_s$ is a projector.

Now consider the left $S^2_{qs}$-module $(S^2_{qs})^2e_s :=
\{(x,y)e_s \; | \; x,y\in
S_{qs}^2\}$, and the linear map
$\Theta_s: (S^2_{qs})^2e_s\to E_s$
given by $\Theta_s: (x,y)e_s\mapsto x(\alpha+s\beta) +y(\gamma+s\delta)$. This
map is well-defined because $(x,y)e_s =0$ if and only if
\begin{equation}
\left\{\begin{array}{ll}
      (x(\alpha+s\beta)+y(\gamma+s\delta))(\delta -qs\gamma) = 0\\
      (x(\alpha+s\beta)+y(\gamma+s\delta))(s\alpha -q^{-1}\beta) = 0
\end{array} \right.
\label{brace}
\end{equation}
Multiplying the first equation by $\alpha+s\beta$ and the second one by
$\gamma+s\delta$ and then using (\ref{fact}) one deduces that
$x(\alpha+s\beta)
+y(\gamma+s\delta)= 0$. The map $\Theta_s$ is clearly surjective.
Also, if $x(\alpha+s\beta)
+y(\gamma+s\delta) =0$, then (\ref{brace}) are
satisfied, so that $(x,y)e_s = 0$. This implies that $\Theta_s$ is an
isomorphism of left $S^2_{qs}$-modules. Therefore, $E_s$ is
isomorphic to $(S^2_{qs})^2e_s$ and hence it
is a finitely generated projective module.

Next we compute the Chern number or the monopole charge of $E_s$.
The formula (4.4) in
\cite{MasNak:non} for the Chern character of $S^2_{qs}$ gives
$\tau^1(1) =0$, $\tau^1(\zeta) = {1+s^2\over 1-q^{-2}}$.
Hence  the Chern number of $E_s$  is ${\rm ch}(E_s) :=
\tau^1({\rm tr} e_s) =
\tau^1(1 - (\zeta -q^{-2}\zeta)/(1+s^2)) =-1$, and is independent of
$s$. By the same argument
as in \cite{HajMaj:pro} we conclude that
$E_s$ is not isomorphic to $S_{qs}^2\otimes V$. This implies
\cite{Brz:mod}
that $SU_q(2)$ is
not isomorphic to $S_{qs}^2\otimes C_s$ as a left
$S_{qs}^2$-module and a right
$C_s$-comodule (i.e. $SU_q(2)(S_{qs}^2)^{C_s}$ is not {\em cleft}).
This means that
the coalgebra Hopf bundle of $S_{qs}^2$ is not trivial.

Finally we show that the {\em monopole connection} in
$SU_q(S^2_{qs})^{C_s}$, constructed in
\cite{BrzMaj:geo}, is the Grassmannian connection in $E_s$. Its
connection one-form is given by
\begin{equation}
\omega(g^\pm_n) = Si(g^\pm_n)\sw 1di(g^\pm_n)\sw 2,
\label{mon}
\end{equation}
where
$$
i(g^+_n)=\prod_{k=0}^{n-1}{\alpha + q^ks(\beta+\gamma) + q^{2k}s^2\delta
\over 1+q^{2k}s^2}, \quad
i(g^-_n)=\prod_{k=0}^{n-1}{\delta - q^{-k}s(\beta+\gamma) +
q^{-2k}s^2\alpha \over 1+q^{-2k}s^2}.
$$
The connection (\ref{mon}) is strong and the corresponding covariant
derivative $\nabla$ on $E_s$ computed from
(\ref{cov.der}) is
$\nabla(u) = du -u\omega(g_1^+)$, for all $u\in E_s$. Let $u =
x(\alpha+s\beta) + y(\gamma +s\delta)$ for some $x,y\in S^2_{qs}$.
Using the explicit form of $\omega(g^+_1)$ we have
\begin{eqnarray*}
\nabla(u) & = & d(x(\alpha+s\beta)) +d(y(\gamma+s\delta))\\
&&\!\!\!\!\!\!\!\!\!-{1\over 1+s^2}(x(\alpha+s\beta) + y(\gamma +s\delta))
((\delta-qs\gamma)d(\alpha+s\beta)
+(s\alpha -q^{-1}\beta)d(\gamma+s\delta)).
\end{eqnarray*}
With the help of the Leibniz rule and (\ref{fact}) this can be gathered
in the following form
\begin{eqnarray*}
\nabla(u) & = & (d(x,y)\! +\!(x,y) {1\over 2}d\pmatrix{(\alpha+s\beta)
(\delta -qs\gamma) &
(\alpha+s\beta)(s\alpha -q^{-1}\beta) \cr (\gamma+s\delta)
(\delta -qs\gamma) &
(\gamma+s\delta)(s\alpha - q^{-1}\beta)})\pmatrix{\alpha \!+ \!s\beta \cr
\gamma \! + \! s\delta}\\
&=& ((dx,dy) +(x,y) de_s)\pmatrix{\alpha+s\beta \cr \gamma +s\delta}.
\end{eqnarray*}
Now viewing $u$ in $(S^2_{qs})^2e_s$ via $\Theta^{-1}_s$ and $\nabla$
as a map $(S^2_{qs})^2e_s\to \Omega^1S^2_{qs}\otimes
_{S^2_{qs}} (S^2_{qs})^2e_s$ via $\Theta_s$ one obtains
\begin{equation}
\nabla((x,y)e_s) = (d(x,y) +(x,y)de_s)e_s.
\label{gras}
\end{equation}
The covariant derivative of the form (\ref{gras}) is called the
{\em Grassmannian connection} on a projective module. Therefore,
similarly to the standard quantum sphere case discussed in
\cite{HajMaj:pro},  the covariant
derivative corresponding to the $q$-monopole connection on any quantum
sphere $S^2_{qs}$
is the Grassmannian connection on $E_s$.

\end{document}